\documentclass[10pt]{article}
\usepackage{amssymb,amsmath,amsthm,latexsym}
\usepackage{amsfonts}
\usepackage[twlrm]{rawfonts}
\usepackage[english]{babel}
\usepackage[latin1]{inputenc}
\usepackage{fancyhdr}
\usepackage{indentfirst}
\usepackage{color}
\newtheorem{defi}{Definition}[section]

\newtheorem{teo}{Theorem}[section]
\newtheorem{cor}{Corolary}[section]
\newtheorem{lem}{Lemma}[section]
\newtheorem{pro}{Proposition}[section]

\newenvironment{dem}[1][Proof]{\noindent\textbf{#1.} }{\hfill \rule{0.5em}{0.5em}}

\newcommand{\N}{\mathbb{N}}

\newcommand{\R}{\mathbb{R}}

\newcommand{\C}{C(\overline{\Omega})}

\newtheorem{ob}{Remark}[section]

\newtheorem{Af}{Claim}[section]

\begin{document}
	
\setlength{\baselineskip}{6.5mm} \setlength{\oddsidemargin}{8mm}
\setlength{\topmargin}{-3mm}

\title{\bf An Ambrosetti-Prodi type result for integral equations involving dispersal operator}

\author {Natan de Assis Lima$^a$\thanks{N. A. Lima, natan.mat@cche.uepb.edu.br}\, , \,    Marco A. S. Souto$^b$\thanks{M. A. S. Souto was partially supported by CNPq/Brazil 306082/2017-9 and INCT-MAT, marco@mat.ufcg.edu.br}\,\,\,\,\,\,\vspace{2mm}
\and {\small $a.$ Universidade Estadual da Paraí­ba} \\ {\small Centro de Ciências Humanas e Exatas} \\ {\small CEP: 58500-000, Monteiro - Pb, Brazil}\\
{\small $b.$ Universidade Federal de Campina Grande} \\ {\small Unidade Acad\^emica de Matem\'{a}tica} \\ {\small CEP: 58429-900, Campina Grande - Pb, Brazil}}

%\author{Claudianor O. Alves\thanks{C. O. Alves was partially supported by CNPq/Brazil xxxxxx/xxxx-x and INCT-MAT, coalves@dme.ufcg.edu.br}\, , Natan de Assis Lima\thanks{N. A. Lima, natan.mat@cche.uepb.edu.br}\, and \,Marco A. S. Souto \thanks{M. A. S. Souto was partially supported by CNPq/Brazil
%		xxxxxx/xxxx-x and INCT-MAT, marco@dme.ufcg.edu.br}\,\,\,\,\,\,\vspace{2mm}
%		\and {\small  Universidade Federal de Campina Grande} \\ {\small Unidade AcadÃªmica de MatemÃ¡tica} \\ {\small CEP: 58429-900, Campina Grande - PB, Brazil}\\}

\date{}

\maketitle

\begin{abstract}
	In this paper we study the existence of solution for the following class of nonlocal problems
\[
L_0u =f(x,u)+g(x) , \ \mbox{in} \ \Omega,
\]
where $\Omega \subset \R^{N}$, $N\geq 1$, is a bounded connected open, $g \in \C$, $f:\overline{\Omega} \times \R \to \R$ are function, and $L_0 : \C \to \C$ is a nonlocal dispersal operator. Using a sub-supersolution method and the degree theory for $\gamma$-Condensing maps, we have obtained a result of the Ambrosetti-Prodi type, that is, we obtain a necessary condition on $g$ for the non-existence of solutions, the existence of at least one solution, and the existence of at least two distinct solutions.


\end{abstract}

\section{Introduction}

In this work we study the existence of solutions for the following equation 
\[
L_0u =f(x,u)+g(x) , \ \mbox{in} \ \Omega, \eqno{(P)}
\]
where $\Omega \subset \R^{N}$, $N\geq 1$, is a smooth bounded domain, $g \in \C$, $f:\overline{\Omega} \times \R \to \R$ is a function that verifies some hypotheses that will be detailed below, and $L_0 : \C \to \C$ is the nonlocal dispersal operator given by 
\begin{equation}\label{eq02}
L_0 u(x) = \int_{\Omega}K(x,y)u(y)dy, \quad \mbox{for} \ u \in \C,
\end{equation} 
with a continuous and nonnegative dispersal kernel $K$. The dispersion mechanism is currently a focus of theoretical interest and recently has received much attention. Most of these continuous dispersion models are based on reaction-diffusion equations, which are widely studied see \cite{Bates-Chmaj}, \cite{Bates-Fife-Ren-Wang}, \cite{Bates-Zhao}, \cite{Chasseigne-Chaves-Rossi}, \cite{Chen}, \cite{Coville 2}, \cite{Coville 3}, \cite{Garcia-Rossi}, \cite{Kao-Lou-Shen}, \cite{Kao-Lou-Shen 2}, \cite{Rossi}, \cite{Sun-Shi-Wang}. This type of diffusion process has been broadly used to describe the dispersion of a population (of cell or organisms) through the environment, as are indicated in \cite{Fife 1}, \cite{Fife 2}, \cite{Hutson-Martinez-Mischaikow-Vickers}, if $u(y)$ is thought of as a density at a location $y$, $K(x,y)$ as the probability distribution of jumping from a location $y$ to a location $x$, then the rate at which the individuals from all other places are arriving at location $x$ is $$\int_{\Omega}K(x,y)u(y)dy.$$ The presence of nonlocal reaction term in equation $(P)$ means, from the biological point of view, that the crowding effect depends not only on their own point in space but also depends on the entire population in an $N$-dimensional habitat $\Omega$, see \cite{Furter-Grinfeld}. In many problems in biology (and ecology), for example in \cite{Berestycki-Coville-Hoang}, Berestycki, Coville and Hoang-Hung are interested in finding persistence criteria for a species that has a long range dispersal strategy. For such a specific model species, we can think of trees of which seeds and pollens are disseminated on a wide range. The possibility of a long range dispersal is well known in ecology, where numerous data now available support this assumptions, see \cite{Cain-Milligan-Strand}, \cite{Clark} and \cite{Schurr-Steinitz-Nathan}. For further dispersion problems, with this dispersion formulation of individuals, see also \cite{Medlock-Kot} and \cite{Murray}.

The motivation to study $(P)$ comes from the Ambrosetti-Prodi result \cite{Ambrosetti-Prodi}, which is a very interesting result concerning the solvability of the following Dirichlet problem 
\begin{equation}\label{eq03}
\left\{
\begin{array}{lcl}
-\Delta u = f(u) + g(x) \quad \mbox{in} \ \Omega, \\
u=0 \quad \mbox{on} \ \partial \Omega,
\end{array}
\right.
\end{equation}
where $\Omega$ is a bounded domain in $\R^N$ with a boundary $\partial \Omega$ of class $C^{2,\alpha}$. The nonlinearity is given by a $C^2$ real-valued function $f: \R \to \R$ such that $f''(s)>0$ for all $s \in \R$ and 
\begin{equation}\label{eq04}
0< \lim_{s\rightarrow -\infty}f'(s) < \lambda_1 < \lim_{s\rightarrow +\infty}f'(s) < \lambda_2
\end{equation}
where $\lambda_1$ and $\lambda_2$ are eigenvalues of $(-\Delta, H^{1}_0)$. Under these assumptions they proved the existence of a $C^1$ manifold $\cal{M}$ in $C^{0,\alpha}(\overline{\Omega})$ which split this space in two open sets $O_0$ and $O_2$ with the following property: $(i)$ if $g \in O_0$, problem (\ref{eq03}) has no solution; $(ii)$ if $g \in \cal{M}$, problem (\ref{eq03}) has exactly one solution; $(iii)$ if $g \in O_2$, problem (\ref{eq03}) has exactly two solutions. In \cite{Ambrosetti-Prodi}, the method that was used is based on inversion theorems for differentiable mappings with singularities in Banach spaces. This method is quite beautiful and geometrical, but it seems to depend heavily on the fact that $f$ is a convex function. Furthermore, the work of Ambrosetti and Prodi has the inconvenient of not giving necessary or sufficient conditions to $(i)$, $(ii)$ and $(iii)$ occur.

In $1975$, Beger and Podolak \cite{Berger-Podolak}, have taken a major step in the study of the problem and obtained a cartesian structure for the $\cal{M}$ manifold in the Hilbert spaces. They decomposed the functions $g \in C^{0,\alpha}(\overline{\Omega})$ in the form $g = t \phi_1 + g_1$, where $\phi_1$ is a normalized positive eigenfunction (in $L^2(\Omega)$) associated with the eigenvalue $\lambda_1$ (eigenvalue of $(-\Delta, H^{1}_0)$) and $g_1 \in \{\phi_1\}^\perp$ in the $L^2(\Omega)$ sense, i.e., $\displaystyle\int_{\Omega}g_1(x)\phi_1(x)dx =0$. Thus, they wrote equation (\ref{eq03}) as
\begin{equation}\label{eq05}
\left\{
\begin{array}{lcl}
-\Delta u = f(u) + t \phi_1(x) + g_1(x) \quad \mbox{in} \ \Omega, \\
u=0 \quad \mbox{on} \ \partial \Omega,
\end{array}
\right.
\end{equation}
Using the Lyapunov-Schimidt method, for each $g_1$ as above, they found a real number $t(g_1) \in \R$ depending continuously on $g_1$, so that: $(i)$ $g \in O_0$ (i.e. (\ref{eq05}) has no solution), if $t > t(g_1)$; $(ii)$ $g \in \cal{M}$, (i.e. (\ref{eq05}) has exactly one solution), if $t=t(g_1)$; $(iii)$ $g \in O_2$, (i.e. (\ref{eq05}) has exactly two solutions), if $t<t(g_1)$.

Also in $1975$, Kazdan and Warner \cite{Kazdan-Warner}, published a long article dealing with uniform second order elliptic operators with Dirichlet or Newmann conditions. They worked by substituting hypotheses (\ref{eq04}) for hypotheses
\begin{equation}\label{eq06}
-\infty\leq \limsup_{s\rightarrow -\infty}\dfrac{f(s)}{s} < \lambda_1 < \liminf_{s\rightarrow +\infty}\dfrac{f(s)}{s} \leq +\infty
\end{equation}
which does not involve the derivative of $f$. They found a sub and a supersolution to $t$ negative enough, and using the method of monotone iteration, proved the existence of a solution. Even removing the convexity of function $f$, proved the existence of a function $t: \{\phi_1\}^\perp \to \R$ such that: $(i)$ (\ref{eq05}) has no solution, if $t > t(g_1)$; $(ii)$ (\ref{eq05}) has at least one solution, if $t<t(g_1)$.

Posteriorly, Amann and Hess \cite{Amann-Hess}, and concomitantly Dancer \cite{Dancer}, improved the work of Kazdan and Warner by finding at least two solutions to $t<t(g_1)$, and at least one solution to $t=t(g_1)$. They have used theory of degree arguments to obtain this result.

It is observed that the strict convexity of the function $f$ implies in the possibility of obtaining in each case the exact number of solutions. On the other hand, the position of the limits $\displaystyle\lim_{|s| \rightarrow +\infty}\frac{f(s)}{s}$ in relation to the spectrum of $(-\Delta, H^1_{0})$ (how many eigenvalues does this limit cut) influences the number of solutions we can obtain. For example, if, in addition to the hypotheses relevant to problem (\ref{eq05}), we assume the convexity (\ref{eq06}) and the hypotheses that $\displaystyle\lim_{s\rightarrow +\infty}\frac{f(s)}{s}< \lambda_2$, then for each $g_1 \in \{\phi_1\}^\perp$ (in the $L^2(\Omega)$ sense), there exists $t(g_1)$ such that: $(i)$ (\ref{eq05}) has no solution, if $t > t(g_1)$; $(ii)$ (\ref{eq05}) has exactly one solution, if $t=t(g_1)$; $(iii)$ (\ref{eq05}) has exactly two solutions), if $t<t(g_1)$. 

The proof of this result is very close to the idea developed by Berestycki \cite{Berestycki} and can be found in Figueiredo \cite{Figueiredo}.

Another example, if we consider the hypotheses
\[
\limsup_{s\rightarrow -\infty}\dfrac{f(s)}{s} < \lambda_1 < \lambda_2 < \lim_{s\rightarrow +\infty}\dfrac{f(s)}{s} <\lambda_3,
\]
then one can find $\tau \in \R$ such that (\ref{eq05}) has at least three solutions if $t<\tau$.

In the present paper, our main objective is to show  the existence of a result of Ambrosetti-Prodi type for the problem $(P)$. That is, Ambrosetti-Prodi type problems are characterized by the determination of the functions $g$ for which the $(P)$ equation has solution or not, and if so, the minimum number (or, if possible, exact number) of solutions.

Observe that the problem $(P)$ can be written as follows 
\[
L_0 u = f(x,u) + t \phi_1 + g_1(x) \quad \mbox{in} \ \Omega,  \eqno{(P)_t}
\]
where $g \in \C$ is decomposed as $g(x)=t\phi_1(x)+g_1(x), x \in \Omega$, where $\phi_1$ is a normalized positive eigenfunction associated to the principal eigenvalue $\lambda_1$ of $L_0$, and $g_1 \in \{\phi_1\}^\perp$, i.e., $\displaystyle\int_{\Omega}g_1(x)\phi_1(x)dx =0$.

In all work the $K: \overline{\Omega} \times \overline{\Omega} \to \R$ is a nonnegative function that verifies
\begin{itemize}
\item [$(K_{1})$] $K(x,y)=K(y,x)$ for all $x,y\in \overline\Omega$;
\item [$(K_{2})$] There exists $\delta>0$ such that $K(x,y) > 0$ for all $x,y \in \overline{\Omega}$ and $|x-y|\leq \delta$.
\end{itemize}

So, this paper is methodically organized as we indicated next: In Section 2, we show a more general version of the Maximum Principle for $L_0$ found in $\cite{Garcia-Rossi}$. In Section 3 supposing that $f:\overline{\Omega} \times \R \to \R$ is a locally Lipschitz function and is increasing with respect to variable $t \in \R$ such that:
\begin{itemize}
\item[$(f_1)$] There exist $A>\|k\|_\infty$ and $C>0$ such that $f(x,s) \geq As-C$ for all $s \geq 0$ and for all $x \in \overline{\Omega}$, where $k(x)=\displaystyle\int_{\Omega}K(x,y)dy$;
\end{itemize}
This function $k$ has an important role in our results. When $k(x)\equiv 1$ we have the case that appears in \cite{Chasseigne-Chaves-Rossi}, \cite{Garcia-Rossi 2} and \cite{Rossi}.

With this hypothesis we find a number $m>0$ in which the problem $(P)_t$ has no positive solution if $t>m$. In addition, if assume that $f$ also verifies
\begin{itemize}
\item[$(f_2)$] There is a number $0<a<\lambda_1$ such that $\displaystyle\lim_{s\rightarrow -\infty}\frac{f(x,s)}{s} = a$, for all $x \in \overline{\Omega}$.
\end{itemize}
we find a number $m>0$ such that problem $(P)_t$ has no solution (positive, negative or that changes the signal) if $t>m$.

Note that all function that we have verified $\displaystyle\liminf_{s\rightarrow \infty}\frac{f(x,s)}{s}>\|k\|_\infty$ also satisfies $(f_1)$.

In Section 4 we study the existence of a solution for equation $(P)_t$ without boundary conditions. Assuming that $ f $ satisfies the condition $(f_1)$ we show the existence of at least one positive solution. More precisely we have the following result: 
\begin{teo} \label{Teo 1}
Assume $(K_1)$ and $(K_2)$, and suppose that $f$ is a locally Lipschitz function, increasing with respect to variable $t \in \R$ e verifies $(f_1)$.
Then for all $g_1 \in \{\phi_1\}^\perp$, there is a real number $t(g_1)$ such that
\begin{itemize}
\item[$(i)$] the problem $(P)_t$ has no positive solution, if $t>t(g_1)$;
\item[$(ii)$] the problem $(P)_t$ has at least one positive solution, if $t\leq t(g_1)$.
\end{itemize}
\end{teo}

Here, the existence assertion is proved by the monotone iteration method.

In Section 5 we have analyzed the existence of a second solution for the problem $(P)$. In order to obtain this second solution, we have assumed in $f$ the condition $(f_2)$ and the condition 
\begin{itemize}
\item[$(f_3)$] For all compact $\cal{K} \subset \R$ there is a number $\sigma>0$ such that $\dfrac{f(x,s)-f(x,t)}{s-t}> \sigma$, for all $s,t \in \cal{K}$ and all $x\in \overline{\Omega}$.
\end{itemize}
Note that all function that verifies $f_{t}'(x,t)>0$, for all $x\in \overline{\Omega}$, also satisfies $(f_3)$.

With the above hypotheses, we have the main result in this paper.

\begin{teo} \label{Teo 2}
Under the hypotheses of Theorem \ref{Teo 1}, in addition to $(f_2)$ and $(f_3)$ we have that, for all $g_1 \in \{\phi_1\}^\perp$, there is a real number $t(g_1)$ such that
\begin{itemize}
\item[$(i)$] the problem $(P)_t$ has at least two solutions, if $t<t(g_1)$;

\item[$(ii)$] the problem $(P)_t$ has at least one solution, if $t=t(g_1)$.
\end{itemize}

\end{teo}

In Section 5, since the difference of the problem with $L_0$ from the problem with $-\Delta$ is that, $(-\Delta)^{-1}$ is a compact operator, while $(L_0-AI)^{-1}$ is not a compact for any $A\geq 0$ large, we can not use the Leray-Schauder degree theory. Faced with this problem we will make use of the degree theory for $\gamma$-Condensing maps (see Deimling \cite{Deimling}, pp.$71$), which is an extension of the Leray-Schauder degree, being a larger class of perturbations of identity, defined in terms of noncompactness measures, since there are several interesting types of functional equations that can not be treated by compact operators, but need the structure of this larger class. 

We have made some comments in the final of the Section 4 that we have at most one solution for $(P)_t$ to any $t\in \R$, when $f$ just satisfies the condition $(f_1)$ and the condition below,
\begin{itemize}
\item[$(f_4)$] $|f(x,s)-f(x,t)| > \|k\|_\infty |s-t|$ for all $s,t \in \R_+$, uniformly in $x \in \overline{\Omega}$.
\end{itemize}

\begin{ob}
An interesting fact is that we can replace $\phi_1$ by the constant $\phi_0 \equiv 1$ and consider the problem
\[
L_0 u = f(x,u) + t + g_1(x) \quad \mbox{in} \ \Omega,  \eqno{(Q)_t}
\]
where $t \in \R$, $g_1 \in \C$ is such that $\displaystyle\int_{\Omega}g_1(x)dx=0$, and $f$ verifies the hypotheses $(f_1)$ and $(f_2)$, we have that the proof of Theorem \ref{Teo 1} and Theorem \ref{Teo 2} follows in an analogous way of the problem $(P)_t$.
\end{ob}

\section{A Maximum Principle}
 
In this section, we consider some preliminary facts related to the $L_0$ operator. Here, we are assuming that $K$ is a nonnegative continuous function that verifies $(K_1)$ and $(K_2)$. 

\begin{ob} \label{rm2}
 	The function  $k:\overline\Omega \to \mathbb R$	given by
 	\[
 	k(x)=\int_\Omega K(x,y)dy,
 	\]
 	is very important to the  our maximum principle. It is easy to check that
 	\[
 	\int_{\Omega\times \Omega}K(x,y)v(y)^2dxdy= \int_\Omega k(x)v(x)^2dx, 
 	\]
 and \[\langle L_0v, v \rangle \leq \|k\|_\infty \|v\|^{2}_2\] where $ \langle u,v\rangle=\displaystyle\int_\Omega uv dx$ is the inner product of $L^2(\Omega)$. Since $K$ is symmetric, we also have
 	\begin{eqnarray*}
 		\int_{\Omega\times \Omega}K(x,y)[v(x)-v(y)]^2dxdy=2\int_\Omega k(x)v(x)^2dx\\ -2\int_{\Omega\times \Omega}K(x,y)v(x)v(y)dxdy 
 	\end{eqnarray*}
 \end{ob}

The next result sets out details about the principal eigenvalue of $L_0$ and it can be found in \cite{Alves-Lima-Souto} (to see also \cite{Garcia-Rossi 2}).

\begin{pro} \label{Prop.1}
	The eigenvalue problem 
	\[
	L_0 u = \lambda u, \quad \mbox{in} \ \Omega, 
	\]
	has an unique eigenvalue $\lambda_1 > 0$ whose the eigenfunction are continuous on $\overline{\Omega}$ with defined signal and $\dim N(L-\lambda_1 I)=1$. Moreover, $\lambda_1 =\sup \sigma(L_0)$, where $\sigma(L_0)$ is spectrum of operator $L_0:\C \to\C$.
\end{pro}

\begin{lem}
Let $\lambda_1$ be the principal eigenvalue of $L_0$. Then, $\lambda_1 \equiv k(x)$ or $\|k\|_\infty > \lambda_1 > \displaystyle\inf_{x \in \overline{\Omega}} k(x)$.
\end{lem}

\begin{dem}
In fact, $\langle L_0\phi_1,1 \rangle = \langle \phi_1, L_0 1\rangle$ by the symmetry of $L_0$, thus
\[
\lambda_1 \int_{\Omega}\phi_1(x)dx = \int_{\Omega}k(x)\phi_1(x)dx
\]
hence,
\[
\int_{\Omega} (\lambda_1 - k(x))\phi_1(x)dx=0.
\]
The proof is done.

\end{dem}

As the Laplacian operator, the dispersal operator also possesses an important maximum principle (here, we are showing a version that generalizes the one found in \cite{Garcia-Rossi}). 

\begin{lem}\label{Pmax}(Maximum Priciple) Suppose that $K$ satisfies $(K_1)$, $(K_2)$ and $\Omega$ is a connected open set.
	Let $u \in C(\overline{\Omega})$ verify $$L_0u(x) - c(x)u(x)\leq 0,$$ for all $x\in\Omega$, where $c $ is a bounded function satisfying $c(x)> k(x)$ a.e. in $\Omega$. Then $u > 0$ or $u \equiv 0$ in $\Omega$.
\end{lem}

\begin{dem}
 Let $u^{-}(x)=\min\{u,0\}$. Note that, $u^{-}(x) \leq 0,$ for all $x \in \Omega$. Since $L_0u(x) - c(x)u(x) \leq 0$ in $\Omega$,
$$\int_{\Omega}\left(\int_{\Omega}K(x,y) u(y)dy - c(x)u(x)\right)u^{-}(x)dx \geq 0,$$
$$\int_{\Omega\times \Omega}K(x,y) u(y)u^{-}(x)dydx \geq \int_{\Omega}c(x)u(x)u^{-}(x)dx,$$
$$\int_{\Omega\times \Omega}K(x,y) u(y)u^{-}(x)dydx \geq \int_{\Omega}c(x)u^{-}(x)^{2}dx \mbox{ and }$$
$$\int_{\Omega\times \Omega}K(x,y) u^{+}(y)u^{-}(x)dydx + \int_{\Omega\times \Omega}K(x,y) u^{-}(y)u^{-}(x)dydx \geq \int_{\Omega}c(x)u^{-}(x)^2dx.$$
Since $$\displaystyle\int_{\Omega\times\Omega}K(x,y) u^{+}(y)u^{-}(x)dydx \leq 0$$ we have
$$\int_{\Omega\times \Omega}K(x,y) u^{-}(y)u^{-}(x)dydx \geq \int_{\Omega}c(x)u^{-}(x)^{2}dx.$$
From Remark \ref{rm2}
$$-\frac 12\int_{\Omega\times\Omega}K(x,y) (u^-(x)-u^-(y))^{2}dydx + \int_{\Omega}k(x)u^-(x)^{2}dx\geq \int_{\Omega}c(x)u^{-}(x)^{2}dx$$
or,
$$0\geq-\frac 12\int_{\Omega\times\Omega}K(x,y) (u^-(x)-u^-(y))^{2}dydx \geq \int_{\Omega}[c(x)-k(x)]u^{-}(x)^{2}dx\geq0.
$$
We must have $u^-\equiv0$ in $\Omega$,  that is, $u \geq 0$  in $\Omega$.

Furthermore, if $u(x_{1})=0$ for some $x_{1} \in \Omega$, then $$\int_{\Omega}K(x_{1},y)u(y)dy = 0$$
which implies that $u \equiv 0$ in a neighborhood of $x_{1}$. A standard connectedness argument gives $u \equiv 0$ in $\Omega$, we conclude the proof.

\end{dem}

\begin{ob} \label{Op. S}
	As a corollary to the Maximum Priciple, if $w \in \C$  verifies $L_{0}w-c(x)w = 0$ in $\Omega$, where $c $ is a bounded function satisfying $c(x)> k(x)$ a.e. in $\Omega$, then $w \equiv 0$. Thus, for every $f \in L^{2}(\Omega)$, the  problem $L_{0}u-c(x)u = f, \mbox{in}\ \Omega$ admits at most one solution $u \in \C$.
\end{ob}

It is well known that the sub and supersolution method has been used widely in reaction-difusion equations with nonlocal term. What is not different for the dispersal operator. For this purpose, we consider the following system
\begin{equation} \label{eq07}
L_0u  = f(x,u)\mbox { in }  \Omega
\end{equation}
where $f: \overline{\Omega}\times \mathbb{R}\longrightarrow \mathbb{R}$ is a locally Lipschitz function.

First we give the definition of sub and supersolution for the above system.

\begin{defi}
A positive function $\overline{u} \in C(\overline{\Omega})$ is said to be a supersolution of $(\ref{eq07})$, if $$L_0\overline{u} \leq f(x,\overline{u}), \ in \ \Omega.$$
A subsolution $\underline{u} \in C(\overline{\Omega})$ is defined similarly by reversing the inequality.

\end{defi}

The next lemma is related to the regularity of the solutions of (\ref{eq07}).

\begin{lem} (Regularity) \label{regularity}
If $f:\overline{\Omega} \times \R \to \R$ is a locally Lipschitz function and is increasing with respect to variable $t \in \R$ and $u \in L^{\infty}(\Omega)$ verifies
\[
L_0u(x)=f(x,u(x)) \quad \mbox{for all } x \in \Omega,
\]
then $u \in \C$.
\end{lem}

\begin{dem}
We know that, $v(x)=L_0u(x)$ is a continuous function in $x \in \overline{\Omega}$, because $L_0$ is a linear and compact operator, with continuous kernel $K$. Fix $x_0 \in \overline{\Omega}$ and consider $(x_n) \subset \overline{\Omega}$ such that $\displaystyle\lim_{n \longrightarrow \infty}x_n =x_0$.

Thus, since $|u(x_n)| \leq \|u\|_\infty$ for all $n \in \N$, from Bolzano-Weierstrass, there is $(u(x_{n_k})) \subset (u(x_n))$ convergent subsequence, let us say that $$u(x_{n_k}) \to s \in [-\|u_n\|_\infty,\|u_n\|_\infty].$$ Hence, as $v(x_{n_k})=L_0u(x_{n_k})=f(x_{n_k},u(x_{n_k}))$ by taking the limit of $k \to \infty$, we get $v(x_0)=f(x_0,s)$. On the other hand, $$f(x_0,u(x_0))=L_0u(x_0)=v(x_0)=f(x_0,s).$$ Therefore, since $f$ is increasing we have $s=u(x_0)$, that is, $u(x_{n_k}) \to u(x_0)$ in $\C$. Hence, $u \in \C$.

\end{dem}

\begin{lem} \label{Ite.Mon.}
	Suppose that $(\ref{eq07})$ has a positive supersolution $\overline{u}$ and a positive subsolution $\underline{u}$ defined on $\Omega$ such that $\underline{u} \leq \overline{u}$. Besides, suppose that $f:\overline{\Omega} \times \mathbb{R}\to \mathbb{R}$ is a locally Lipschitz function and is increasing with respect to variable $t \in \R$. Then $(\ref{eq07})$ has a solution $u \in \C$ satisfying $\underline{u} \leq u \leq \overline{u}$.
	
\end{lem}

\begin{dem}
Define $\Sigma = \{u \in \C; \underline{u}\leq u \leq \overline{u} \}$. From Lemma \ref{Pmax}, we have that $L_0v(x)-\beta v(x)$ admits, at most, a single solution, if $\beta > k(x)$ a.e. in $\Omega$. Moreover, since $t \longmapsto f(x,t)$ is increasing with respect to variable $t \in \R$, we have that $t \longmapsto f(x,t)-\beta t$ is decreasing in $t \in \R$, for each fixed $x \in \overline{\Omega}$.

Now, from Remark \ref{Op. S}, for each $u \in \C$ the problem
\[
L_0v(x)-\beta v(x)=f(x,u(x))-\beta u(x), \quad \mbox{for all } x \in \Omega
\]
admits at most one solution $v \in \C$. 

Next, if $w_1,w_2 \in \Sigma$ are such that $w_1 \leq w_2$,
$$
L_0v_1(x)-\beta v_1(x)=f(x,w_1(x))-\beta w_1(x), \quad \mbox{for all } x \in \Omega
$$
and
$$
L_0v_2(x)-\beta v_2(x)=f(x,w_2(x))-\beta w_2(x), \quad \mbox{for all } x \in \Omega,
$$
then $$f(x,w_2(x))-\beta w_2(x) \leq f(x,w_1(x))-\beta w_1(x), \quad \mbox{for all } x \in \Omega.$$ This implies that, $$L_0v_2(x)-\beta v_2(x) \leq L_0v_1(x)-\beta v_1(x), \quad \mbox{for all } x \in \Omega,$$ that is, $$L_0(v_2(x)-v_1(x))-\beta(v_2(x)-v_1(x)) \leq 0, \quad \mbox{for all } x \in \Omega$$ from maximum principle (Lemma \ref{Pmax}), we have $v_1(x) \leq v_2(x)$ for all $x \in \Omega$.
   
Consider the sequence $\{u_{n}\}^{\infty}_{n=1} \subset \Sigma$ given by 
$$
 L_0u_1(x)-\beta u_1(x)=f(x,\overline{u}(x))-\beta \overline{u}(x), \quad \mbox{for all } x \in \Omega
$$
$$
 L_0u_2(x)-\beta u_2(x)=f(x,u_1(x))-\beta u_1(x), \quad \mbox{for all } x \in \Omega
$$
and
$$
 L_0u_n(x)-\beta u_n(x)=f(x,u_{n-1}(x))-\beta u_{n-1}(x), \quad \mbox{for all } x \in \Omega
$$
which is a monotone decreasing sequence. Similarly, we also obtain another sequence $\{v_{n}\}^{\infty}_{n=1} \subset \Sigma$ by 
$$
L_0v_1(x)-\beta v_1(x)=f(x,\underline{u}(x))-\beta \underline{u}(x), \quad \mbox{for all } x \in \Omega
$$
$$
L_0v_2(x)-\beta v_2(x)=f(x,v_1(x))-\beta v_1(x), \quad \mbox{for all } x \in \Omega
$$
and
$$
L_0v_n(x)-\beta v_n(x)=f(x,v_{n-1}(x))-\beta v_{n-1}(x), \quad \mbox{for all } x \in \Omega
$$
which is a monotone increasing sequence. Furthermore, we have 
    $$
    \underline{u}\leq v_{1} \leq ...\leq v_{n}\leq u_{n} \leq ... \leq u_{1} \leq \overline{u},
    $$
     then, there exist functions $u^{*}, v^{*}$ such that 
     $$
     u^{*}(x) = \lim_{n \rightarrow \infty} u_{n}(x) \ \text{and} \ v^{*}(x) = \lim_{n \rightarrow \infty} v_{n}(x) \quad \mbox{pointwisely in } \Omega.
     $$
     It follows that, $\underline{u} \leq v^{*} \leq u^{*} \leq \overline{u}$. By the Lebesgue's Dominated Convergence Theorem, we have $$\lim_{n\rightarrow \infty} L_0 u_{n}(x) = L_0 u^{*}(x),$$ $$\lim_{n\rightarrow \infty} L_0 v_{n}(x) = L_0 v^{*}(x).$$ 
On the other hand, $(L_0u_n)$ and $(L_0v_n)$ are uniformly convergent in $\C$, we assume that $L_0u_n \to w$ and $L_0v_n \to z$ in $\C$ respectively. As $L_0$ is a linear and compact operator, we have 
$$
L_0u_n(x) = \int_{\Omega}K(x,y)u_n(y)dy \rightarrow \int_{\Omega}K(x,y)u^*(y)dy = L_0u^*(x) \quad \mbox{in}\ \Omega,
$$
and
$$
L_0v_n(x) = \int_{\Omega}K(x,y)v_n(y)dy \rightarrow \int_{\Omega}K(x,y)v^*(y)dy = L_0v^*(x) \quad \mbox{in}\ \Omega,
$$
so, $L_0u_n \to L_0u^*$  in $\C$ and $L_0v_n \to L_0v^*$  in $\C$.
     
Due to the continuity of $f$, $$\lim_{n\rightarrow \infty} [f(x,u_{n}(x)) + \beta u_{n}(x)] = f(x,u^{*}(x)) + \beta u^{*}(x),$$ $$\lim_{n\rightarrow \infty} [f(x,v_{n}(x)) + \beta v_{n}(x)] = f(x,v^{*}(x)) + \beta v^{*}(x).$$

Hence, we have $$L_0 u^{*}(x) - \beta u^{*}(x) = f(x,u^{*}) - \beta u^{*}(x), \quad \mbox{ for all } x \in \Omega$$ $$L_0 v^{*}(x) - \beta v^{*}(x) = f(x,v^{*})- \beta v^{*}(x), \quad \mbox{ for all } x \in \Omega$$ that is, $u^{*}$ and $v^{*}$ are solutions of the problem (\ref{eq07}). Therefore, from Lemma \ref{regularity}, we have that $u^*, v^*\in \C$ and the proof is done.

\end{dem}
 
\section{Nonnexistence of solution }

In this section we are going to prove the nonnexistence of solution for large $t \in \R$. Consider problem $(P)_t$, that is
\[
L_0u = f(x,u) +t\phi_1 + g_1(x), \quad x \in \Omega.
\]
We see that $u$ solves $(P)$ if, and only if, $u$ solves $(P)_t$.

First, let us show the nonexistence of positive solution.

\begin{lem} \label{mp}
Assume $(f_1)$. Then there exists $m>0$, independent of $g_1 \in \{\phi_1\}^\perp$ such that, for all $t>m$, problem $(P)_t$ has no positive solution.
\end{lem}

\begin{dem}
Suppose that problem $(P)_t$ has a positive solution $u$, thus $L_0u=f(x,u)+t\phi_1 +g_1$ in $\Omega$. Multiplying the above equation by $\phi_1$ and integrating into $\Omega$ we obtain $$\int_{\Omega}\phi_1L_0udx = \int_{\Omega}\phi_1 f(x,u)dx + t\int_{\Omega}\phi_1^{2}dx +\int_{\Omega}\phi_1g_1dx$$ thus $$\int_{\Omega}\phi_1L_0udx = \int_{\Omega}\phi_1 f(x,u)dx + t$$ from $(f_1)$ we have $$\int_{\Omega}\phi_1L_0udx \geq A\int_{\Omega}\phi_1 udx - C\int_{\Omega}\phi_1dx+ t$$ since $L_0$ is symmetric, $$\lambda_1\int_{\Omega}\phi_1udx \geq A\int_{\Omega}\phi_1 udx - C\int_{\Omega}\phi_1dx+ t.$$ Hence, $$t \leq (\lambda_1- A)\int_{\Omega}\phi_1 udx + C\int_{\Omega}\phi_1dx \leq C\int_{\Omega}\phi_1dx.$$
Therefore, the existence of positive solution $u$ to $(P)_t$ necessarily implies that $$t\leq m:=C\int_{\Omega}\phi_1dx,$$ that is, if $t>m$ problem $(P)_t$ has no solution.

\end{dem}

This result does not affirm that, for $t>m$, the problem $(P)_t$ has no negative or nodal solution. Also does not affirm that, for $t\leq m$, the problem $(P)_t$ has positive solution.

Now, assuming $f$ also satisfies $(f_2)$ we still can find a $m>0$ in which the problem $(P)_t$ has no solution (positive, negative or nodal) to $t>m$. Before we need the following estimate,

\begin{ob}
Assume $(f_1)$ and $(f_2)$. Then there is a number $C_1>0$ such that for all $x \in \overline{\Omega}$
\begin{equation} \label{Desig. f}
f(s) \geq As-C_1, \quad \mbox{for all } s \in \R \quad \mbox{and}\quad f(s) \geq (a+\epsilon)s-C_1, \quad \mbox{for all } s \in \R.
\end{equation}
\end{ob}

\begin{lem} \label{m}
Assume $(f_1)$ and $(f_2)$. Then there exist $m>0$, independent of $g_1 \in \{\phi_1\}^\perp$ such that, for all $t>m$, problem $(P)_t$ has no solution (positive, negative or nodal).
\end{lem}

\begin{dem} Similarly to Lemma \ref{mp}, using (\ref{Desig. f}) we have
$$t \leq (\lambda_1- A)\int_{\Omega}\phi_1 udx + C_1\int_{\Omega}\phi_1dx \leq C_1\int_{\Omega}\phi_1dx, \quad \mbox{if } \int_{\Omega}u\phi_1dx \geq 0$$ and $$t \leq (\lambda_1 - (a+\epsilon))\int_{\Omega}\phi_1 udx + C_1\int_{\Omega}\phi_1dx \leq C_1\int_{\Omega}\phi_1dx, \quad \mbox{if } \int_{\Omega}u\phi_1dx \leq 0.$$
Therefore, the existence of solution $u$ to $(P)_t$ necessarily implies that $$t\leq m:=C_1\int_{\Omega}\phi_1dx,$$ that is, if $t>m$ problem $(P)_t$ has no solution.

\end{dem}

\section{Proof of Theorem \ref{Teo 1}}

In this section, we will show the existence of a positive solution to the problem $(P)_t$, that is, prove the Theorem \ref{Teo 1}. Here, we are assuming that $f$ only verifies $(f_1)$. 

The next lemma will guarantee the existence of a supersolution for $(P)$.

\begin{lem} \label{Super}
Assume $(f_1)$. Then, for all $g \in \C$, problem $(P)$ has a supersolution $w \in \C$. Furthermore, any subsolution $u \in \C$ of $(P)$ is such that $u<w$ in $\overline{\Omega}$.
\end{lem}

\begin{dem}
Let $w \in \C$ be the unique solution of
\[
L_0w - Aw = -\|g\|_\infty - C, \mbox{ in } \Omega,
\]
where $A$ and $C$ have been introduced in property $(f_1)$. Thus,
\[
L_0w = Aw -\|g\|_\infty -C < f(x,w) + g(x), \mbox{ in } \Omega
\]
that is, $w \in \C$ is a supersolution of $(P)$. Since $A> \|k\|_\infty$ from Lemma \ref{Pmax}, we have $w>0$ in $\overline{\Omega}$.

Now, let $u \in \C$ be a subsolution of $(P)$, that is
\[
L_0u \geq f(x,u) +g(x), \mbox{ in } \Omega.
\]
Therefore, $$L_0(w-u)< Aw-\|g\|_\infty -C - f(x,u) -g(x)$$ or, $$L_0(w-u)< A(w-u)-\|g\|_\infty -g(x)<A(w-u)$$which implies $$L_0(w-u)-A(w-u)<0.$$ Hence, from Lemma \ref{Pmax}, we have that $w>u$ in $\overline{\Omega}$.

\end{dem}

\begin{cor} \label{Cor.Super}
Assume $(f_1)$. Let $g_1 \in \{\phi_1\}^\perp$ and $t \in \R$ be given, then there exists $R>0$ such that, if $u \in \C$ is a positive function such that $$L_0u(x)=f(x,u)+t\phi_1(x)+g_1(x), \quad \mbox{for all } x \in \Omega,$$ we have that $\|u\|_\infty<R$.
\end{cor}

\begin{dem}
From Lemma \ref{Super}, there exists $w \in \C$ supersolution of $(P)_t$ such that $0<u < w$ in $\overline{\Omega}$. Hence, \[\|u\|_\infty < \|w\|_\infty := R.\] Therefore, there exists $R>0$ such that $\|u\|_\infty < R$.
\end{dem}

\begin{ob} \label{Ob.Super}
If $(P)$ has solution for $g \in \C$, then for all $h \in \C$ such that $h\leq g$ in $\Omega$, problem 
\begin{equation} \label{eq08}
L_0v=f(x,v)+h(x), \mbox{ in }\Omega
\end{equation}
also has solution.

Indeed, let $v \in \C$ be a solution of $(P)$ for given $g\in \C$. Then, $v$ is a subsolution of $(\ref{eq08})$ because \[L_0v - f(x,v)=g(x)\geq h(x), \quad x \in \Omega.\] From Lemma $\ref{Super}$, $(\ref{eq08})$ has a supersolution $w \in \C$ such that, $v<w$ in $\overline{\Omega}$. Hence, from Lemma $\ref{Ite.Mon.}$, $(\ref{eq08})$ has a solution $U\in \C$ such that $v \leq U \leq w$.
\end{ob}

With this observation, the following lemma is proved:

\begin{lem} \label{mm}
Let $g_1 \in \{\phi_1\}^\perp$ be given. Suppose that problem $(P)_t$ has a solution for a given $t_0 \in \R$. Then the problem $(P)_t$ has a solution for any $t<t_0$.
\end{lem}

In order to prove Theorem \ref{Teo 1}, we still need to obtain a subsolution for $(P)_t$. The next lemma will guarantee the existence of this subsolution.

\begin{lem} \label{Sub}
Under the condition $(f_1)$ and given $g_1 \in \{\phi_1\}^\perp$, there exists $t \in \R$ such that problem $(P)_t$ has a subsolution.
\end{lem}

\begin{dem}
Choose $t$ so that $-f(x,0) > t\phi_1(x)+g_1(x)$ for all $x \in \Omega$. Consider $z = \epsilon \phi_1$, thus $$L_0(\epsilon \phi_1) - f(x,\epsilon \phi_1)-t\phi_1-g_1 = \epsilon\lambda_1-f(x,\epsilon \phi_1)-t\phi_1-g_1> 0,$$ for $\epsilon\approx 0^+$ and for all $x \in \Omega$, that is $$L_0z > f(x,z) +t\phi_1+g_1, \quad \mbox{in } \Omega$$ concluding that $z$ is a positive subsolution of $(P)_t$.

\end{dem}

\noindent {\bf Proof of Theorem \ref{Teo 1}}

Let $g_1 \in \{\phi_1\}^\perp$ be given. From Lemma \ref{Sub}, there is $t \in \R$ such that problem $L_0u=f(x,u)+t\phi_1+g_1(x)$ for all $x \in \Omega$ has a positive subsolution $z \in \C$. On the other hand, for these $g_1$ and $t$, from Lemma \ref{Super}, we have a positive supersolution $w \in \C$ such that $z \leq w$ in $\overline{\Omega}$. Thus, from Lemma \ref{Ite.Mon.}, $(P)_t$ has a positive solution $u \in \C$ such that $z \leq u \leq w$.

Thereby, the set $$\Sigma = \{t \in \R; (P)_t \mbox{ has positive solution}\}$$ is not empty, from Lemma \ref{mp} is bounded from above and from Lemma \ref{mm} this set is seen to be a half-line. That is, taking $t(g_1)$ to be the supremum of the $\Sigma$, it follows that for all $t>t(g_1)$ problem $(P)_t$ has no solution. 

To conclude, we need to show for $t=t(g_1)$. Consider a sequence $t_n< t(g_1)$ such that $t_n \to t(g_1)$. It follows from the above information that the problem $(P)_t$ has a solution $u_n \in \C$ for the given function $g_1$ and each $t_n$, that is
\begin{equation} \label{fn}
L_0u_n = f(x,u_n)+t_n\phi_1 +g_1, \mbox{ in } \Omega.
\end{equation}
Since $(u_n)$ is bounded in $\C$, then $(u_n)$ is bounded in $L^2(\Omega)$. That is, there is a some subsequence of $(u_n)$, still denoted by itself, such that $u_n \rightharpoonup u$ in $L^2(\Omega)$. As $(L_0u_n)$ is uniformly convergent in $\C$, we assume that $L_0u_n \to w$. As $L_0$ is a linear and compact operator, we have 
$$
L_0u_n(x) = \int_{\Omega}K(x,y)u_n(y)dy \rightarrow \int_{\Omega}K(x,y)u(y)dy = L_0u(x) \quad \mbox{in}\ \Omega,
$$ 
and so, $L_0u_n \to L_0u$  in $\C$. On the other hand, from (\ref{fn}) we obtain \[f(x,u_n) = L_0u_n -t_n\phi_1 -g_1 \to L_0u - t_0\phi_1-g_1 := z \mbox{ uniformly in } \Omega. \] From Corolary \ref{Cor.Super}, there exists $R>0$ such that $\|u_n\|_\infty <R,$ for all $n \in \N$. Now, since $f$ is a increasing function with respect to variable $t \in \R$, it follows that there is $\sigma>0$ such that, if $u_n(x)\neq u_m(x)$ we have
\[
\sigma< \inf_{-R\leq u_n(x), u_m(x)\leq R}\dfrac{|f(x,u_n(x))-f(x,u_m(x))|}{|u_n(x)-u_m(x)|}.
\] 
Hence, \[\sigma \|u_n - u_m\|_\infty < \left(\dfrac{|f(x,u_n(x))-f(x,u_m(x))|}{|u_n(x)-u_m(x)|}\|u_n-u_m\|_\infty \right)\leq \|f(.,u_n) - f(.,u_m)\|_\infty \] which implies \[\sigma \|u_n - u_m\|_\infty < \|f(.,u_n) - f(.,u_m)\|_\infty \leq \|f(.,u_n)-z\|_\infty + \|z-f(.,u_m)\|_\infty\]

Therefore, $(u_n)$ is a Cauchy sequence in $\C$. Hence, $$f(.,u_n) \to f(.,u) \mbox{ in } \C.$$ Clearly $u$ is a solution of problem $(P)_t$ with the given $g_1$ and $t=t(g_1)$. Therefore the proof of Theorem \ref{Teo 1} is complete.

\textbf{Comments:} Considering $f$ on the hypotheses $(f_1)$ and $(f_4)$, we can show the uniqueness of solution for each $t \in \R$.

Indeed, consider $u,w \in \C$ positive functions such that
\[
L_0u(x)=f(x,u)+t\phi_1(x)+g_1(x), \quad \mbox{for all } x \in \Omega
\]
and
\[
L_0w(x)=f(x,w)+t\phi_1(x)+g_1(x), \quad \mbox{for all } x \in \Omega.
\]
Hence, if $u(x)\neq w(x)$ we have
\[
L_0(u(x)-w(x)= f(x,u(x))-f(x,w(x))=\left(\dfrac{f(x,u(x))-f(x,w(x))}{u(x)-w(x)}\right)(u(x)-w(x)),
\]
thus
\[
L_0(u(x)-w(x))-a(x)(u(x)-w(x))=0,
\]
where 
\[
a(x) = \left\{
\begin{array}{lcl}
\dfrac{f(x,u(x))-f(x,w(x))}{u(x)-w(x)}, \ \mbox{if } u(x)\neq w(x)  \\
0, \ \mbox{if } \ u(x)=w(x),
\end{array}
\right.
\]
Moreover, from $(f_4)$ we have $a(x)>k(x)$ for all $x \in \overline{\Omega}$. From Lemma \ref{Pmax}, we must have $u\equiv w$ in $\overline{\Omega}$.

\section{Proof of Theorem \ref{Teo 2}}

In this section, in order to obtain a second solution to problem $(P)$, is necessary to make other assumptions under $f$, that is, assuming that $f$ verifies $(f_1)$, $(f_2)$ and $(f_3)$. Moreover, to prove Theorem \ref{Teo 2} we are going to use the theory of degree for $\gamma$-Condensing maps that extend the Leray-Schauder degree to a larger class of perturbations of identities, defined in terms of noncompactness of measures, (see Deimling \cite{Deimling}, pp.$71$). 

We begin with the theorem that establishes a priori estimate:

\begin{teo} (A priori estimate) \label{Priori}
Given $g \in \C$, there is a number $R>0$ such that, if $u$ is a solution of $(P)$, that is, $u$ is a solution of $$L_0 u = f(x,u) + g(x), \quad \mbox{in } \Omega$$ then $\|u\|_\infty<R$.
\end{teo}

\begin{dem}
We know that from $(f_1)$, there is a function $w \in \C$ such that all solution of $(P)$ verifies $u<w$ in $\overline{\Omega}$. Suppose that there is $(u_n) \subset \C$ such that,
\[
\|u_n\|_\infty \longrightarrow \infty \mbox{ and } L_0u_n=f(x,u_n)+g(x), \ \mbox{ in } \Omega.
\]
We have that $u_n <w$ in $\overline{\Omega}$, for all $n \in \N$. Consider $v_n = \dfrac{u_n}{\|u_n\|_\infty}$, then $(v_n)$ is a bounded sequence in $\C \subset L^2(\Omega)$, without loss of generality, we can suppose there is $v \in L^{2}(\Omega)$ such that  $v_n \rightharpoonup v \ \mbox{in} \ L^2(\Omega)$. As $(L_0v_n)$ is uniformly convergent in $\C$, we assume that $L_0v_n \to z$ in $\C$. But since $L_0$ is a linear and compact operator, we have 
$$
L_0v_n(x) = \int_{\Omega}K(x,y)v_n(y)dy \rightarrow \int_{\Omega}K(x,y)v(y)dy = L_0v(x) \quad \mbox{in}\ \Omega,
$$ 
and so, $L_0v_n \to L_0v$  in $\C$.

Now, from $(f_2)$ for all $\epsilon>0$ there is $C_\epsilon>0$ such that for all $x \in \overline{\Omega}$
$$(a-\epsilon)s+C_\epsilon \geq f(x,s)\geq (a+\epsilon)s-C_\epsilon,\quad \mbox{for all }  s \leq \|w\|_\infty.$$ Thus, $(a-\epsilon)u_n+g(x)+C_\epsilon \geq f(x,s)\geq (a+\epsilon)u_n+g(x)-C_\epsilon$ which implies
\[
(a-\epsilon)v_n+\dfrac{g(x)+C_\epsilon}{\|u_n\|_\infty} \geq f(x,s)\geq (a+\epsilon)v_n+\dfrac{g(x)-C_\epsilon}{\|u_n\|_\infty}.
\]
Since $\|u_n\|_\infty \longrightarrow \infty$ and $u_n <w$ for all $n \in \N$, then there is a $x_n \in \overline{\Omega}$ such that $$-u_n(x_n) = |u_n(x_n)| = \|u_n\|_\infty.$$ Hence,
\[
(a-\epsilon)=(a-\epsilon)\|v_n\|_\infty = (a-\epsilon)|v_n(x_n)| =-(a-\epsilon)v_n(x_n)\leq -L_0v_n(x_n)+\dfrac{g(x_n)+C_\epsilon}{\|u_n\|_\infty}
\]
that is, $$(a-\epsilon) \leq \|L_0v_n\|_\infty + \dfrac{\|g\|_\infty + C_\epsilon}{\|u_n\|_\infty}.$$ Choosing $\epsilon = \dfrac{a}{2}$ and passing the limit, we see that $\|L_0v\|_\infty \geq \dfrac{a}{2}$. Therefore, $\|v\|_\infty>0$, $v \neq 0$ and $L_0v=av$ in $\Omega$. On the other hand, we have
$$
v_n(x)=\dfrac{u_n(x)}{\|u_n\|_\infty} < \dfrac{w(x)}{\|u_n\|_\infty} \rightarrow 0, \quad \mbox{when } n\rightarrow \infty, \ \mbox{for all } x \in \Omega,
$$
which implies $v(x)=\lim_{n \rightarrow \infty}v_n(x)<0$ for all $x \in \Omega$. Thus, we have $v(x)<0$ for all $x \in \Omega$, that is, $v$ would have defined signal and would be autofunction associated with the eigenvalue $a < \lambda_1$, which is absurd.

\end{dem}

Now, for each $t \in \R$ we define an operator $F_t : \C \to \C$ given by
\begin{equation} \label{P2}
F_t u:= \dfrac{1}{M}L_0u +u - \dfrac{1}{M}f(x,u)- \dfrac{1}{M}(t\phi_1 + g_1(x)), \mbox{ for some } M>0.
\end{equation}
Note that, fixed point for this operator is solution to problem $(P)_t$. Indeed, if $u \in \C$ is such that $F_tu=u$ we have
\[
u = F_tu= \dfrac{1}{M}L_0u +u - \dfrac{1}{M}f(x,u)- \dfrac{1}{M}(t\phi_1 + g_1(x))
\]
or
\[
L_0u = f(x,u)+ t\phi_1 + g_1(x), \ \mbox{in } \Omega
\]
that is, $u$ is a solution of $(P)_t$.

With the above information, we have defined a continuous mapping $F_t$, that is not compact. In addition, we have previously that the solutions of our problem $(P)_t$ are precisely the zeros of $I-F_t$. To obtain our result, we need to show that the $F_t$ function is $\gamma$-Condensing.

\begin{ob}
From Theorem \ref{Priori}, there is a $R>0$ such that, if $u$ is a solution of $(P)_t$ then $\|u\|_\infty<R$. On the other hand, since $f$ is a locally Lipschitz function, there is a number $M>0$ such that $$M>\Gamma =\sup_{-R\leq s,t\leq R}\left(\dfrac{f(x,s)-f(x,t)}{s-t}\right), \ \mbox{ for all } x \in \overline{\Omega}.$$ From $(f_3)$, there is $\sigma>0$ such that
\[
0<\sigma<\dfrac{f(x,s)-f(x,t)}{s-t}\leq \Gamma, \mbox{ for all } -R \leq s,t \leq R,  \ \mbox{ for all } x \in \overline{\Omega}.
\]
With the above study we obtain
\begin{equation} \label{f}
0<\left(1-\dfrac{f(x,s)-f(x,t)}{M(s-t)}\right)<1-\dfrac{\sigma}{M}, \mbox{ for all } -R\leq s,t\leq R,  \ \mbox{ for all } x \in \overline{\Omega}.
\end{equation}

\end{ob}

\begin{defi}
Let $X$ be a Banach space and $\cal{B}$ it is bounded sets. Then $\alpha : \cal{B} \to \R_+$, defined by
\[
\alpha(B)=\inf \{d>0; B \ \mbox{admits a finite cover by sets of diameter } \leq d \},
\]
is called the \textit{Kuratowski-measure of noncompactness}, and $\beta : \cal{B} \to \R_+$ defined by
\[
\beta(B)=\inf \{r>0; B \ \mbox{can be covered by finitely many balls of radius } r \},
\]
is called the \textit{ball measure of noncompactness}.
\end{defi}

\begin{ob}
If $\gamma: \cal{B} \to \R_+$ be either $\alpha$ or $\beta$, then
\begin{itemize}
\item [$(a)$] $\gamma(B)=0$ iff $\overline{B}$ is compact;
\item[$(b)$] $\gamma$ is a seminorm, i.e., $\gamma(\lambda B)=\lambda \gamma(B)$ and $\gamma(B_1+B_2) \leq \gamma(B_1)+\gamma(B_2)$
\end{itemize}
The proof of these items can be found in Deimling \cite{Deimling}, pp.$41$.
\end{ob}

\begin{lem}
Supose that $(f_1)$, $(f_2)$ and $(f_3)$ hold. Then the operator $F_t$ is an application $\gamma$-Condensing, for all $ t \in \R$.
\end{lem}

\begin{dem}
We say that an application is $\gamma$-condensing when $\gamma(F(B))<\gamma(B)$ whenever $B \subset \C$ is bounded and $\gamma(B)>0$, where $\gamma$ will be the \textit{Kuratowski-measure} $\alpha$ or the \textit{ball noncompactness measure} $\beta$. Thus, if $B \subset \C$ is bounded, we have $\gamma(F_t(B)) \leq \gamma(L_0(B)) + \gamma(G(x,B))$ where $G(x,s) = s-\dfrac{1}{M}f(x,s)$ in $x \in \overline{\Omega}$. But, since $L_0$ is compact we must have $\gamma(L_0(B))=0$. On the other hand, $G$ is a contraction for all $s \in I \subset \R$ and all $x \in \overline{\Omega}$. Indeed, note that, for $-R\leq u(x),v(x)\leq R$ (for some $R>0$) we have from (\ref{f}) and for all $x \in \overline{\Omega}$, \[0< \left|1-\dfrac{f(x,u(x))-f(x,v(x))}{M(u(x)-v(x)}\right| < 1-\dfrac{\sigma}{M}, \mbox{ for } u(x) \neq v(x).\]
Hence, \[|G(x,u(x))-G(x,v(x))| = |u(x)-f(x,u(x))-v(x)+f(x,v(x))|, \quad x \in \overline{\Omega}\] or \[|G(x,u(x))-G(x,v(x))| \leq \left|1-\dfrac{f(x,u(x))-f(x,v(x))}{M(u(x)-v(x)}\right| |u(x)-v(x)|, \quad x \in \overline{\Omega} \] which implies \[|G(x,u(x))-G(x,v(x))| < (1-\dfrac{\sigma}{M})|u(x)-v(x)|, \quad x \in \overline{\Omega}. \] For this last inequality, we are using (\ref{f}). Therefore, $G$ is contraction for all $x \in \overline{\Omega}$, so $\gamma(G(x,B)) < \gamma(B)$, for $B=B_R(0) \subset \C$ bounded.

With the study made above we have to $\gamma(F_t(B))<\gamma(B)$ where $B=B_R(0) \subset \C$, i.e., $F_t$ is an application $\gamma$-Condensing, for all $t \in \R$.

\end{dem}

Now we determine the degree of $I-F_t$ for a certain subset of $\C$.

\begin{lem} \label{BR}
Suppose that $(f_1)$, $(f_2)$ and $(f_3)$ hold. Let $g_1 \in \{\phi_1\}^\perp$ and $t_0 < t(g_1)$ be given. Then there is an $R>0$ such that \[d(I-F_{t_0}, B_R, 0) = 0,  \] where $B_R= \{u \in \C; \|u\|_\infty < R \}$.
\end{lem}

\begin{dem}
From Theorem \ref{Teo 1} the problem $(P)_t$ has no solution if $t>t(g_1)$. Thus, choose a $t_1>t(g_1)$. It follows from the Theorem \ref{Priori} that there is a constant $R>0$ such that $\|u\|_\infty < R$ for all eventual solution of problem $(P)_{t_0}$ with $g_1$ fixed. It is also follows that this inequality occur for all eventual solution of problem $(P)_t$ for all $t \in [t_0, t_1]$, from Lemma \ref{mm}. Since $F_t$, with $t \in [t_0, t_1]$, constitutes an admissible homotopy in between $F_{t_0}$ and $F_{t_1}$, because \[(I-F_t)u \neq 0, \mbox{ for all } \|u\|_\infty = R \mbox{ and all } t \in [t_0, t_1], \] we have $d(I-F_{t_0}, B_R, 0) = d(I-F_{t_1}, B_R, 0)$. But $d(I-F_{t_1}, B_R, 0) = 0$ because problem $(P)_t$ has no solution for $t=t_1$. Therefore, $d(I-F_{t_0}, B_R, 0) = 0$.

\end{dem}

\begin{lem} \label{Open}
Suppose that $(f_1)$, $(f_2)$ and $(f_3)$ hold. Let $g_1 \in \{\phi_1\}^\perp$ and $t_0 < t(g_1)$ be given. Then there exist $M>0$ and an open $U \subset \C$ such that \[d(I-F_{t_0}, U, 0) = 1.\]
\end{lem}

\begin{dem}
It follows from there exists $v \in \C$ which is positive solution of $(P)_{t_1}$, with $t_0 < t_1< t(g_1)$. Moreover, $v$ is a subsolution of $(P)_t$ when $t=t_0$, that is
\[L_0v - f(x,v) = t_1\phi_1 +g_1 > t_0\phi_1 +g_1, \mbox{ in } \Omega  \] or 
\begin{equation} \label{IV}
L_0v > f(x,v) + t_0\phi_1 +g_1, \mbox{ in } \Omega
\end{equation}
On the other hand, from Lemma \ref{Super}, there is $w \in \C$ supersolution of $(P)_{t_0}$, that is
\begin{equation} \label{V}
L_0w < f(x,w) + t_0\phi_1 +g_1, \mbox{ in } \Omega
\end{equation}
Moreover, $v < w$ in $\overline{\Omega}$.

Now, choose $M>0$ so that, $M>\|k\|_\infty$, $f(x,s)-Ms$ is a decreasing function in $0 \leq s \leq \|w\|_\infty$ and that
\[
F_{t_0} u:= \dfrac{1}{M}L_0u +u - \dfrac{1}{M}f(x,u)- \dfrac{1}{M}(t\phi_1 + g_1(x)), \ \mbox{ for } x \in \overline{\Omega}
\]
is $\gamma$-Condensing map.

Define $W = \{u \in \C; v<u<w \mbox{ in } \overline{\Omega} \}$, we have that $W$ is open, bounded and convex in $\C$.

\begin{Af}
$F_{t_0}: \overline{W} \to \C$ is such that $F_{t_0}(\overline{W}) \subset W$. 
\end{Af}
In fact, if $u \in \overline{W}$ then $v\leq u \leq w$ in $\overline{\Omega}$. Let $z=F_t(u)$, thus $$Mz=L_0u+(Mu-f(x,u))-t_0\phi_1-g_1, \mbox{ in } \Omega.$$ Now, note that
\[
L_0v - Mv > f(x,v) -Mv + t_0\phi_1 +g_1, \mbox{ in } \Omega
\]
which implies
\[
Mv < L_0v + (Mv -f(x,w))- t_0\phi_1 -g_1, \mbox{ in } \Omega.
\]
In an analogous way, $Mw > L_0w + (Mw -f(x,w))- t_0\phi_1 -g_1$ in $\Omega$. Hence,
\[
M(z-v) > L_0(z-v) + [(Mu-f(x,u))-(Mv-f(x,v))], \mbox{ in } \Omega.
\]
Since $v\leq u$ we have $L_0v \leq L_0u$ and $(Mv-f(x,v)) \leq (Mu-f(x,u))$, thus $M(z-v)>0$ which implies $v <z$ in $\overline{\Omega}$. Similarly it is proved that $z<w$ in $\overline{\Omega}$. Therefore, $z \in W$.

With the above, we conclude that, if $u \in \partial W$ then $u \neq F_{t_0}(u)$, because if $u \in \overline{W}$ and $u = F_{t_0}(u)$ we must have $u \in W$.

By the study done, $d(I-F_{t_0}, W,0)$ is well defined. Now let us calculate your value. Consider $\psi= \dfrac{v+w}{2}$ we have $v<\psi<w$, i.e. $\psi \in W$. Define $H_\theta(u) = (1-\theta)F_{t_0}(u)+\theta \psi$. For $0\leq \theta \leq 1$ we have $H_\theta : \overline{W} \to W$.

Indeed, we know that if $u \in \overline{W}$ then $F_{t_0}(u) \in W$. Thus, $v<F_{t_0}(u)<w$ and as $v<\psi<w$ we have $H_\theta(\overline{W}) \subset W$, that is $H_\theta$ is a admissible homotopy for all $0 \leq \theta \leq 1$.

Since $u \neq H_\theta(u)$ for all $u \in \partial W$ and all $\theta \in [0,1]$, we conclude that \[ d(I-H_0, W, 0) = D(I-H_1, W, 0). \] But $H_\theta$ is independent of $\theta \in [0,1]$ and  $d(I-H_1, W, 0) = 1$ because that $\psi \in W$. Hence  $d(I-H_0, W, 0) = 1$, and the lemma is proved.

\end{dem}

\noindent {\bf Proof of Theorem \ref{Teo 2}}

$(i)$ Let $g_1 \in \{\phi_1\}^\perp$ and $t_0<t(g_1)$ be given. From Lemma \ref{Open}, there is an $M>0$ and a bounded open set $W$ such that $d(I-F_{t_0}, W, 0) =1$. Thus, $I-F_{t_0}$ has a zero in $W$, that is, problem $(P)_{t_0}$ has a solution $u_1 \in W$, with these given $g_1$ and $t_0$. Now choose $R>0$ such that $W \subset B_R(0)$. From Lemma \ref{BR}, it follows that $d(I-F_{t_0}, B_R,0)=0$, thus $$d(I-F_{t_0}, B_R\setminus \overline{W}, 0) =-1.$$ Hence, problem $(P)_{t_0}$ has another solution $u_2 \in B_R\setminus \overline{W}$. Moreover, $u_1 \neq u_2$.

$(ii)$ Consider a sequence $t_n< t(g_1)$ such that $t_n \to t(g_1)$. It follows from Theorem \ref{Teo 1} that problem $(P)_t$ has a solution $u_n \in \C$ for the given function $g_1$ and each $t_n$, that is
\begin{equation} \label{VII}
L_0u_n = f(x,u_n)+t_n\phi_1 +g_1, \mbox{ in } \Omega.
\end{equation}
Since $(u_n)$ is bounded in $\C$, then $(u_n)$ is bounded in $L^p(\Omega)$, $p>0$. If $p>1$, there is a some subsequence of $(u_n)$, still denoted by itself, such that $u_n \rightharpoonup u$ in $L^p(\Omega)$. As $(L_0u_n)$ is uniformly convergent in $\C$, we assume that $L_0u_n \to w$. As $L_0$ is a linear and compact operator, we have 
$$
L_0u_n(x) = \int_{\Omega}K(x,y)u_n(y)dy \rightarrow \int_{\Omega}K(x,y)u(y)dy = L_0u(x) \quad \mbox{in}\ \Omega,
$$ 
and so, $L_0u_n \to L_0u$  in $\C$. On the other hand, from (\ref{VII}) we obtain \[f(x,u_n) = L_0u_n -t_n\phi_1 -g_1 \to L_0u - t_0\phi_1-g_1 := z \mbox{ uniformly in } \Omega. \] From $(f_3)$ it follows that \[\sigma \|u_n - u_m\|_\infty < |f(x,u_n(x)) - f(x,u_m(x))| \leq \|f(.,u_n)-z\|_\infty + \|z-f(.,u_m)\|_\infty \] Therefore, $(u_n)$ is a Cauchy sequence in $\C$. Hence, $$f(.,u_n) \to f(.,u) \mbox{ in } \C.$$ Clearly $u$ is a solution of problem $(P)_t$ with the given $g_1$ and $t=t(g_1)$. In this context, the proof of Theorem \ref{Teo 2} is complete.

\end{document}